\documentclass[9pt,amstex]{article}
\usepackage{amssymb,amsmath} 
\usepackage{latexsym}

\newtheorem{rmk}{Remark}[section]

\def\buildrel#1_#2^#3{\mathrel{\mathop{\kern 0pt#1}\limits_{#2}^{#3}}}

\newcommand{\Ad}{\mbox{$\mathtt{Ad}$}}
\newcommand{\taub}{{\overline{\tau}}{}} 

\newcommand{\A}{\mathbb A}
\newcommand{\T}{\mathbb T}
\newcommand{\At}{\A_{\theta}}
\newcommand{\C}{\mathbb C} 

\newcommand{\R}{\mathbb R}

\newcommand{\Z}{\mathbb Z} 

\newcommand{\EE}{{\cal E}}

\newcommand{\g}{{\mathfrak{g}}{}}

\newcommand{\n}{{\mathfrak{n}}{}}

\newcommand{\so}{{\mathfrak{so}}{}}

\renewcommand{\a}{{\mathfrak{a}}{}} 
 
\newcommand{\E}{\bf E{}}

\newcommand{\ads}{{\mbox{AdS}_3}{}} 
\newcommand{\sld}{{\mathfrak{sl}(2,\R)}{}} 
 
\newcommand{\OO}{{\cal O}{}}
\newcommand{\HH}{{\cal H}{}}
\renewcommand{\SS}{{\cal S}{}}
\newcommand{\MM}{{\cal M}{}}
\newcommand{\JJ}{{\cal J}{}}
\newcommand{\RR}{{\mathfrak{R}}{}} 
\newcommand{\UU}{{\cal U}{}}

\def\cref#1{Corollary~\ref{#1}}

\title{}

\author{}

\title{Noncommutative Locally Anti-de Sitter Black Holes}
 \author{{\bf P.~Bieliavsky} \\{\it D\'epartement de Math\'ematique}\\ {\it
         Universit\'e Catholique de Louvain}\\ {\it
         Chemin du cyclotron, 2, }\\
         {\it 1348 Louvain-La-Neuve, Belgium}\\ E-mail:
        bieliavsky@math.ucl.ac.be\\
        {\bf S.~Detournay} \thanks{"Chercheur FRIA", Belgium}~ and {\bf Ph.~Spindel}
         \\{\it M\'ecanique et Gravitation}\\ {\it Universit\'e de Mons-Hainaut, 20
         Place du Parc}\\ {\it 7000 Mons, Belgium}\\E-mail:
        stephane.detournay@umh.ac.be,spindel@umh.ac.be\\
        {\bf M.~Rooman} \thanks{FNRS Research Director}\\{\it Service de Physique
         th\'eorique}\\ {\it Universit\'e Libre de Bruxelles, Campus Plaine, C.P.225}\\
         {\it Boulevard du Triomphe, B-1050 Bruxelles, Belgium}\\E-mail:
   mrooman@ulb.ac.be}
\date{January 4, 2005} 

\begin{document}
\maketitle
\begin{abstract}
This note is based on a talk given by one of us at the workshop `Noncommutative Geometry and Physics 2004' (Feb. 2004, Keio University, Japan).
We give a review of our joint work on strict deformation of BHTZ 2+1 black holes.
\end{abstract}
This note is based on a talk given by one of us at the workshop `Noncommutative Geometry and Physics 2004' (Feb. 2004, Keio University, Japan).
We give a review of our joint work on strict deformation of BHTZ 2+1 black holes \cite{BRS02,BDHRS03}.
However some results presented here are not published elsewhere, and an effort is made
for enlightening the instrinsical aspect of the constructions. This shows in particular that
the three dimensional case treated here could be generalized to an  anti-de Sitter space of arbitrary
dimension provided one disposes of a universal deformation formula for the actions of a parabolic 
subgroup of its isometry group. 

\section{Motivations}
The universal covering space of a generic BHTZ space-time \cite{BHTZ93}---realized 
as some open domain in $\ads$---
 is, in a canonical way, the total space of a principal fibration 
over $\R$ with structure group a minimal parabolic
 subgroup $\RR$ of $G:=\widetilde{SL(2,\R)}$ \cite{BRS02}. The action of 
the structure group is isometric w.r.t. the $\ads$ metric on the 
 total space. For spinless BHTZ black holes, the fibers are dense 
open subsets 
of twisted conjugacy classes in $G$ associated to an exterior 
automorphism of $G$. 
These twisted conjugacy classes are known to be WZW branes in 
$G=\ads$ \cite{BP01}
i.e. extremal for the DBI brane action associated to a specific 
2-form $B$ on $\ads$ (referred hereafter as the `$B$-field').
In the present work, we are interested in studying deformations of 
BHTZ spaces which are supported on (or tangential to) the
fibration fibers. We will require our deformations to be non-formal 
(`strict' in the sense of Rieffel, see below) as well as compatible
with the action of the structure group. Our motivation for this is 
threefold.
First, a deformation of the brane in the direction of the $B$-field 
is generally understood as the 
(non-commutative) geometrical framework for studying interactions of 
strings with endpoints attached to the brane \cite{SW99}.
Despite the fact that curved situations have been extensively studied
in various non-commutative contexts, 
in the framework of strict deformation theory i.e. in a purely operator 
algebraic framework, 
non-commutative spaces emerging from string theory have mainly been 
studied in
the case of constant $B$-fields in flat (Minkowski) backgrounds. 
Second, the maximal fibration preserving isometry group of (the 
universal cover of) a BHTZ space being the above mentioned
minimal parabolic group $\RR$, it is natural to ask for a deformation 
which is invariant under 
the action of $\RR$. But there is also a deeper reason for requiring 
the invariance. From considerations
on black hole entropy, or just by geometrical interest, one 
may be interested in defining higher genus locally $\ads$ black holes
(every BHTZ black hole is topologically $S^1\times\R^2$). This would 
of course imply implementing
the action of a Fuchsian group in the BHTZ picture. At the classical 
level, this question was at the
center of important works by Brill et al. \cite{ABBHP98}. Our point here is 
that one may also attend this question at the deformed level.
Indeed, the symmetry group $\RR$ is certainly too small to contain 
large Fuchsian groups. However,
if the deformation is already $\RR$-invariant one may ask wether the 
(classical) action of $\RR$ would 
extend to a (deformed) action of the entire $G$ by automorphisms of 
the deformed algebra. If this is the case (and it is),
one obtains an action of every Fuchsian group at the deformed 
functional level \cite{B99,BP04}.
The third motivation relies on Connes-Landi's and Connes-Dubois-Violette's work 
\cite{CL01,CDV02} where they use 
Rieffel's strict defomation method for actions
of tori to define spectral triples for non-commutative spherical 
manifolds. 
The main point of their construction is that the data of an {\sl 
isometric} action
of a torus on a spin manifold yields not only a strict deformation of 
its function
algebra but a compatible deformation of the Dirac operator as well.
Therefore, disposing of a strict deformation formula for actions of 
$\RR$ 
would produce, in the present situation of BHTZ spaces, examples of 
spectral triples for non-commutative
non-compact Lorentzian manifolds with constant curvature with the 
additional feature that the 
deformation would be supported on the $\RR$-orbits. The difficulty 
here is of course that
these orbits cannot be obtained as orbits of an isometric action of 
$\R^d$.

\section{Generic BHTZ spaces revisited}
A BHTZ black hole is defined as  the quotient space of (an open 
subset of) $\ads$ by an isometric action
of the integers $\Z$. Namely, if $\Xi$ denotes a Killing vector field 
on $\ads$ and if $\{\phi^\Xi_t\}_{t\in\R}$
denotes its flow, the $\Z$-action will be of the form:
\begin{equation}\label{ZACTION}
\Z\times\ads\longrightarrow\ads:(n,x)\mapsto n.x:=\phi^\Xi_n(x).
\end{equation}
Not all Killing fields $\Xi$ give rise to black hole solutions. The 
relevant ones can be defined as follows.
As a Lorentzian manifold, the space $\ads$ is identified with the 
universal cover, $G$, of the group $SL(2,\R)$ of $2\times2$ 
real matrices  of determinant equal to $1$ endowed with its Killing 
metric. The Lie algebra of the isometry 
group of $\ads$ (i.e. the algebra of Killing fields) is then 
canonically isomorphic to the semisimple Lie algebra $\g\oplus\g$ 
where $\g:=\sld$ denotes the Lie algebra of $G$. An element 
$\Xi\in\g\oplus\g$ defines a generic (i.e. non-extremal) BHTZ 
black hole if and only if its adjoint 
orbit only intercepts split Cartan subalgebras of $\g\oplus\g$. The 
squares of the (Killing) norms of the projections of $\Xi$ 
onto each ideals of $\g\oplus\g$ are therefore both positive. Their 
sum
---the Killing norm of $\Xi$--- is called the {\sl mass} and is 
denoted by $M$. Their difference
---which measures how far $\Xi$ is from being conjugated to a 
diagonal element in $\g\oplus\g$---
is called the {\sl angular momentum} and is denoted by $J$.
Generic spaces, as opposed to extremal spaces, correspond to non-zero 
values of $M$ and $J$.
In the present letter, we will mainly be concerned with the two 
following geometrical properties
of generic BHTZ spaces:
\begin{enumerate}
\item[(i)] In a generic BHTZ space, the horizons consist in a finite 
union of  (projected) lateral classes 
of minimal parabolic subgroups of $G$.
\item[(ii)] Every generic BHTZ space is canonically endowed with a 
regular Poisson structure
whose characteristic distribution is (locally) generated by Killing 
vector fields.
\end{enumerate}
We now give a detailed description of the relevant geometry in the 
spinless as well
as in the rotating cases.
\subsection{Spinless BHTZ spaces ($J=0\neq M$)}\label{J=0}
The global geometry in this case can be derived from the two 
following observations.
First, the action (\ref{ZACTION}) appears as the restriction to a 
fixed split (connected) Cartan subgroup $A$
of $G$ of the action of $G$ on itself by twisted conjugation:
\begin{equation}
\tau:G\times G\longrightarrow G:(g,x)\mapsto 
gx\sigma(g^{-1})=:\tau_g(x),
\end{equation}

where $\sigma$ denotes the unique involutive exterior automorphism of 
$G$ 
fixing pointwise the Cartan subgroup $A$. Indeed, if $\a$ denotes the 
Lie subalgebra of $A$ and if 
$H\in\a$ is a vector of Killing length equal to $1$, one has:
\begin{equation}
\phi^\Xi_n=\tau_{\exp(n\sqrt{M}H)}\qquad n\in\Z.
\end{equation}
Second, the action $\tau$ yields the following global decomposition
of $G$. The map
\begin{equation}\label{TIWA}
\phi:A\times G/A\longrightarrow G:(a,gA)\mapsto \tau_g(a)=:\phi(a,gA)
\end{equation}
is well-defined as a global diffeomorphism \footnote{We refer to the 
above decomposition
(\ref{TIWA}) as the {\sl $\sigma$-twisted Iwasawa} decomposition of 
$G$ with repsect to $A$.
Indeed, if $G=KNA$ denotes an Iwasawa decomposition of $G$ associated 
to
the Cartan subgroup $A$, the homogeneous space $G/A$ can be 
identified with
the submanifold $KN$. In this setting, the decomposition (\ref{TIWA}) 
then reads:
\begin{equation*}
K\times N\times A\longrightarrow G:(k,n,a)\mapsto kna\sigma(kn)^{-1},
\end{equation*}
which motivates our terminology.
}.

As a consequence, the space $G=\ads$ appears as the total space of a 
trivial fibration over $A=SO(1,1)\simeq\R$
whose fibers are the $\tau_G$-orbits i.e. the $\sigma$-twisted 
conjugacy classes\footnote{Remark that the foliation of 
$G=\ads$ in $\tau_G$-orbits (\ref{TIWA}) is the unique (up to 
conjugation) foliation tangent to $\Xi$
and generated by Killing vector fields.}. 
As a homogeneous $G$-space every  fiber is isomorphic to the affine 
symmetric space $G/A=\mbox{AdS}_2$. 
Moreover, the $\Z$-action (\ref{ZACTION}) which will define the BHTZ 
space
is fiberwise \footnote{Indeed, through the twisted Iwasawa 
decomposition (\ref{TIWA}) it reads:
\begin{equation}
n.(a,gA)=(a,\exp(n\sqrt{M}H)gA)\qquad n\in\Z.
\end{equation}}. The Killing metric on $G=\ads$ turns out to be 
globally diagonal with respect to  the twisted Iwasawa
decomposition:
\begin{equation}\label{METRIC}
\mbox{d}s^2_{\small{\ads}}=\mbox{d}a^2_{\tiny{A}}\,-\,\frac{1}{4}\cosh^2(a)\,\mbox{d}s^2_{\tiny{G/A}},
\end{equation}

\noindent where $\mbox{d}s^2_{G/A}$ denotes the canonical projected 
$\mbox{AdS}_2$-metric on $G/A$.
The study of the quotient space $\Z\backslash G$ therefore reduces to 
the study of $\Z\backslash (G/A)$.

For a pictural visualization, we realize the space $G/A$ as  the 
$G$-equivariant universal covering space of the adjoint 
orbit $\OO:=\Ad(G)(H)$ of the element $H\in\g$. Sitting in the 
Minkowski space $\g$ endowed with its Killing
form $B$, the orbit $\OO$
is the sphere of radius $1$---a one sheet hyperboloid. The orbits of 
$A$ in $\OO$---along which the $\Z$-identifications will 
be made--- are planar: every intersection 
$\OO\cap(\a^\perp+\xi)\quad\xi\in\OO$ is constituted
by a union of respectively two or five $A$-orbits. The two planar 
sections $\OO\cap(\a^\perp\pm H)=:\SS_\pm$ divide $\OO-\SS_\pm$ 
into six connected open components. In only two of them, $\MM_\pm$, 
the $A$-orbits are timelike curves with respect to the metric 
on $\OO$ induced by the Killing form $B$.
The preimage $\MM$ of $\MM_+\cup\MM_-$ in $A\times G/A=G$ is 
constituted by a countable union of connected open components. 
Denoting by $\MM^x$ the component of $\MM$ containing $x\in\MM$, one 
has
\begin{equation}
\MM=\coprod_{z\in Z(G)}\MM^{zJ}\quad\mbox{ (disjoint union)},
\end{equation}

\noindent where $Z(G)$ denotes the center of $G$ and where  $J$ is 
determined (up to sign) by 
\begin{equation}
J\in K,\qquad J^2=-I.
\end{equation}

\noindent Due to the minus sign in the expression of the Killing 
metric (\ref{METRIC}),  the open set $\MM\subset\ads$ is 
constituted by all the points $x\in\ads$ where the 
Killing vector $\Xi_x$ is spacelike. The $\Z$-action (\ref{ZACTION}) 
is proper and totally discontinuous on $\MM$, 
one therefore has  that the quotient
$\MM^x\to\Z\backslash\MM^x$ defines a metric covering of 
(non-complete) connected Lorentzian manifolds.
These quotient spaces are all isometric to one another and 
homeomorphic to $S^1\times\R^2$. In this picture, the black hole 
singularity appears at the level of $G$ as the preimage $\SS$ of 
$\SS_+\cup\SS_-$ in $A\times G/A$ i.e. the boundary of $\MM$ 
(on $\SS$ the Killing field $\Xi$ is null while it is timelike on the 
complement of the closure 
of $\MM$, therefore yielding closed timelike curves in the quotient). 

At the level of $G$, the horizons may be characterized as follows: a 
point $x\in\MM$ belongs to an horizon if and only 
if the null directions from $x$ which intersect the singularity $\SS$ 
constitute a discrete set among 
all null directions from $x$. The union $\cal H$ of future and past 
horizons are then given by the union of the lateral
classes through the element $J$ of both minimal parabolic subgroups 
of $G$ associated to $A$.
 Formally, if $N$ and $\overline{N}$ denote the nilpotent subgroups 
of $G$ normalized by $A$, one has 
 \begin{equation}
{\cal H}=(J\,Z(G)AN)\,\bigcup\,(J\,Z(G)A\overline{N}),
\end{equation}
\noindent while the singularity is the union of the subgroups:
\begin{equation}
{\cal S}=Z(G)AN\,\bigcup\,Z(G)A\overline{N}.
\end{equation}
 Note that since the conjugation by $J$ defines the $K$-adapted 
Cartan involution,
 the union ${\cal H}$ does not depend on the choice of right or left 
classes. At the 
 level of the black hole (e.g. $ \Z\backslash\MM^J$), the horizons 
are obtained by projecting 
 ${\cal H}$.
 \subsection{Rotating BHTZ black holes ($J\neq0\neq M$)}
 Like in the spinless case, we look for a foliation of $G$ which is 
tangent 
 to $\Xi$ and generated by Killing vector fields. Here again there is 
no choice
 (up to conjugation) for such a foliation and it appears to underly 
an Iwasawa 
 type decomposition of $G$. We adopt the same notations as in the 
 preceding section and we consider $\alpha\in\a^\star$ be such that 
$|\alpha|<1$. 
 Without loss of 
 generality, one may then write the $\Z$-action (\ref{ZACTION}) 
defining the rotating BHTZ space as 
 \begin{equation}
\Z\times G\longrightarrow G:(n,x)\mapsto 
n.x:=\exp(nH)x\exp(n<\alpha,H>\,H).
\end{equation}
 We now observe that the map
 \begin{equation}\label{MICS}
A\times N\times K\longrightarrow G:(a,n,k)\mapsto anka^\alpha
\end{equation}
 is a global diffeomorphism \footnote{for $a\in A$, we set 
$a^\alpha:=\exp(<\alpha,\log(a)>\,\log(a))$.}.

 We call the decomposition (\ref{MICS}) a {\sl modified Iwasawa} 
decomposition. 
 Note that this defines a trivial fibration of $G$ onto $K$ whose 
fibers are
 the orbits of the following action $\taub$ of the solvable group 
$\RR:=AN$:
 \begin{equation}
\taub: \RR\times G\longrightarrow G:(an,x)\mapsto anxa^\alpha.
\end{equation}

Consequently,  every $\taub_\RR$-orbit is stable under the 
$\Z$-action (\ref{ZACTION}) and 
through the modified Iwasawa decomposition (\ref{MICS}) it simply 
reads:
\begin{equation}
n.(r,k)=\taub_{\exp(nH)}(r,k)=(\exp(nH)r,k)\qquad r\in\RR.
\end{equation}
There are three main differences with the $J=0$ case. First, the 
action $\taub$ of $\RR$ does not extend to an action 
of $G$ by isometries. Second, there is no global transversal which is 
everywhere orthogonal to the foliation. The 
$\Z$-action is everywhere proper and discontinuous so that the 
quotient space $\Z\backslash\ads$ is
globally a complete Lorentzian manifold.

The same techniques as in the spinless case  nevertheless apply to 
obtain the causal
structure in the rotating case. Again, one finds that the horizons 
are (projected) lateral 
classes of minimal parabolic subgroups.

\vspace{1cm}

We summarize this section by observing that to every generic BHTZ 
space 
are canonically associated two objects. First, a trivial fibration of 
$\ads$ whose fibers are 
orbits under an isometric 
action of a subgroup of $G$ of dimension at least two. We call it 
{\sl the action associated to
the BHTZ space}. Second, a pair of conjugated minimal parabolic 
subgroups whose 
lateral classes define the horizons. We call them {\sl the parabolic 
subgroups associated 
to the BHTZ space}.

\section{BHTZ domains and extensions}
We consider a generic BHTZ space with corresponding $\Z$-action 
(\ref{ZACTION}) on $\ads$.
We define a {\sl BHTZ domain} of $\ads$ as a $\Z$-stable connected 
simply connected open 
subset $\UU$ of $\ads$ maximal for the following two properties.
\begin{enumerate}
\item[(i)] the restricted action $\Z\times\UU\to\UU$ is proper and 
totally discontinuous;
\item[(ii)] the Killing vector field $\Xi$ is everywhere spacelike on 
$\UU$.
\end{enumerate}

Roughly speaking, a BHTZ domain is thus a realization in $\ads$ of 
the universal covering space of a generic
BHTZ space-time.

\noindent Similarly, we define a {\sl BHTZ extension domain} as a 
$\Z$-stable connected simply connected open 
subset $\EE$of $\ads$ maximal only for property (i) above. We call 
the corresponding quotient manifold
$\Z\backslash\EE$ a {\sl BHTZ extension} of the BHTZ space-time at 
hand. The latter therefore lies in every 
of its maximal extensions.

\noindent Remark that for rotating generic BHTZ spaces, the whole 
$\ads$ constitutes an extension domain.
While for spinless BHTZ spaces, we have seen that any extension 
domain never entirely contains
a $\tau_G$-orbit.

\noindent The rest of this section will be devoted to the following 
assertion:
every BHTZ domain admits a unique (up to conjugation) extension 
domain foliated
by the orbits of the restriction of its associated action to the 
neutral component
of one of its associated minimal parabolic subgroup.

\noindent The above statement is obvious in the rotating case. While 
for the spinless case,
it is just worth observing that at the level of the adjoint orbit 
$\OO\subset\g$ (cf. Subsection \ref{J=0}), the subgroup $AN$ 
(resp. $A\overline{N}$) has exactly 
two open orbits: the connected components of the complement of the 
planar intersection $\OO\cap\n^\perp$ (resp.
$\overline{\n}^\perp$) where $\n$ (resp. $\overline{\n}$) denotes the 
Lie algebra of $N$
(resp. $\overline{N}$). Each of them being identified with the 
subgroup itself (the stabilizer groups
are all trivial), the restricted action of $\exp(\Z\sqrt{M}H)$ is 
proper and discontinuous. One then readily 
verifies that any connected component of the preimage in $A\times 
G/A=G$ of one of these open orbits
defines a BHTZ extension domain.

\section{Poisson structures}
We have seen that the orbits for the associated action of a 
non-spinning
BHTZ space are the $\sigma$-twisted conjugacy classes in $G=\ads$
($\sigma$ is an involutive exterior automorphism of $G$). This remark 
naturally 
leads to consider a particular Poisson structure on $\ads$ whose 
symplectic 
leaves are the $\tau_G$-orbits.
Indeed, on the one hand, these
submanifolds are known to be WZW 1- branes. That is,  extremals of 
the DBI
action associated with the data of a 2-form $B$ on $\ads$ such that 
its exterior
derivative $dB$ is, up to a constant factor, equal to the metric 
volume $\nu$
on $\ads$. On the second hand, each orbit is, as a 
$G$-homogeneous space,
isomorphic to the symmetric space $G/A$ which is endowed with 
a unique (up to a constant factor) $G$-invariant symplectic structure 
$\omega^{G/A}$  \cite{B98}. It is then 
natural to ask wether the restriction of the $B$-field to the tangent 
distribution
of the orbit foliation yields on each orbit (a multiple of) its 
canonical $G$-invariant symplectic structure.
It actually does it, and in a canonical way. Indeed, the twisted 
Iwasawa decomposition
(\ref{TIWA}) $\phi:A\times G/A\to G$ yields a global transversal 
($A$) to the orbit foliation.
By extending by $0$ on $A$ the 2-form $\omega^{G/A}$ to $A\times G/A$ 
one gets
a closed 2-form on $G=\ads$: 
$\omega:=(\phi^{-1})^\star(0\oplus\omega^{G/A})$.
Every $\tau$-invariant 2-form on $G$ is then of the form $f\omega$ 
where $f$
is leafwise constant i.e. $f=f(a)\quad a\in A$. Such a form therefore 
defines an 
admissible $B$-field if and only if $\nu=f'\,da\wedge\omega$. This 
last condition
determines a unique (up to an additive constant) function $f$. 
Indeed, both
$\nu$ and $\omega$ being $\tau$-invariant, it is enough to check the 
condition
on the transversal, which defines $f'$\footnote{ a computation at the 
level of the fundamental vector fields yields
\begin{equation*}
f(a)=\tanh(\frac{a}{2})\,+\mbox{const}.
\end{equation*}
}.

Now for a given generic BHTZ space, the above remark leads us  to 
consider the particular 
class of Poisson structures on $\ads$
constituted by leafwise symplectic structures which are invariant 
under the associated
action. Since the group $\RR$ has a unique (up to a constant factor) 
left-invariant
symplectic structure, such Poisson structures exist in the rotating 
case as well.

\section{Strict deformations for solvable Lie group actions}
\subsection{Connes spectral triples for Abelian Lie group actions}
In this subsection we very roughly recall Rieffel strict deformation 
theory for actions of $\R^d$ \cite{R93} and how
Connes-Landi and Connes-Dubois-Violette used it to define non-commutative spaces 
\cite{CL01,CDV02}.
\subsubsection{Rieffel's machinery}
Given an action $\tau$ of the Abelian group $V:=\R^d$ on a manifold 
${\cal X}$ Rieffel's machinery produces 
a one parameter deformation $\{\At\}_{\theta\in\R}$ of the 
commutative algebra $\A_0:=C({\cal X})$
of complex valued functions on $X$ vanishing at infinity \footnote{if 
$X$ is compact, $C({\cal X})$ is 
the whole algebra of continuous functions.}. The algebras $\At$'s 
involved in the deformation 
are in fact $C^\star$-algebras \footnote{the $C^\star$-norm on $C({\cal X})$ 
is the sup-norm.} 
and the family $\{\At\}_{\theta\in\R}$ is a continuous field of 
$C^\star$-algebras.
One starts with a fixed anti-symmetric matrix $\JJ\in\so(\R^d)$. Note 
that, via the action $\tau$, it determines
a Poisson structure on $C^\infty(X)$: 
$\{u,v\}:=\JJ^{ij}X^\star_i.u\,X^\star_j.v$ where
$u,v\in C^\infty({\cal X})$ and where $\{X^\star_i\}$ are the fundamental 
vector fields w.r.t. the action $\tau$
associated to a basis $\{X_i\}$ of (the Lie algebra of ) $\R^d$. 
Denote by $\alpha:V\times C({\cal X})\to C({\cal X})$
the action on functions induced by $\tau$ and by $\A_0^\infty$  the 
space of smooth vectors of $\alpha$.
One then defines the deformed product at the level of the smooth 
vectors by the following 
oscillatory integral formula:
\begin{equation}\label{RUDF}
a\star_\theta b:=\int_{V\times 
V}e^{ix.y}\alpha_x(a)\alpha_{\theta\JJ{y}}(b)\,\mbox{d}x\mbox{d}y
\end{equation}

where $a,b\in\A_0^\infty$ and where $x.y$ denotes the Euclidean dot 
product on $V=\R^d$ .
The pair $(A_0^\infty,\star_\theta)$ then turns out to be a 
pre-$C^\star$-algebra for a suitable
choice of involution and norm. Its $C^\star$-completion is denoted by 
$\A_\theta$.
Note that for $\theta=0$, the RHS of Formula (\ref{RUDF}) reduces to 
the usual commutative 
product of functions $ab$. Also, the first order expansion term in 
$\theta$ of the oscillatory integral
(\ref{RUDF}) is  $\frac{\theta}{2i}\{a,b\}$. The product 
$a\star_\theta b$ is thus an associative
deformation of $ab$ in the direction of the Poisson bracket coming 
from the action $\tau$.

An important example obtained from Rieffel's machinery is the 
so-called {\sl quantum torus} $C(\T^d)_\theta$.
It is defined as the deformation of the $d$-torus $\T^d$ obtained 
from the above procedure when applied 
to the natural action of $\R^d$ on $\T^d$.

\subsubsection{Deformed spectral triples for compact spin 
$\T^d$-manifolds}
Assume $M$ is a compact spin manifold which the torus $\T^d$ acts 
on by isometries. The point is that Rieffel's machinery applied to 
$M$ only produces a deformation of the {\sl topological} space $M$. 
One 
would want to deform the metric structure as well. For doing this,
Connes introduced the notion of spectral triple. The typical 
commutative example of a spectral triple is the one associated 
to a compact spin manifold; it roughly goes as follows.  
The basic idea is that the metric aspect of the Riemannian manifold  
$M$
is encoded, in the operator algebra framework, by the Dirac operator 
$D$\cite{C94}.
The spectral triple associated to $M$ is then defined as the triple 
$(\A_0,\HH,D)$
where $\A_0=C(M)$ and where $\HH$ is the Hilbert space of 
$L^2$-spinor fields on which 
the Dirac operator $D$ acts as well as the $C^\star$-algebra $\A_0$ 
(simply by multiplication).
Generally, a (not necessarily commutative) spectral triple is defined 
as a triple $(\A,\HH,D)$
where $\A$ is a $C^\star$-algebra, $\HH$ is a Hilbert space 
representation of $\A$ by 
bounded operators and $D$ is a self-adjoint unbounded operator on 
$\HH$ such that,
among other things, $[D,a]$ is bounded for all $a$ in 
$\A$\footnote{this is the essential 
property needed to define the metric distance in operator terms.}.
To produce a deformed spectral triple from an action of $\T^d$ on 
$M$, one first observes that
the deformation $C^\infty(M)_\theta$ obtained from Rieffel's 
machinery can equivalently be  described as
the closed subalgebra $(C^\infty(M)\otimes 
C(\T^d)_\theta)^{\tau\times L^{-1}}$ 
of $C^\infty(M)\otimes C(\T^d)_\theta$
constituted by the elements invariant under the natural action of 
$\T^d$ on the latter tensor 
product algebra ($L$ denotes the regular action). Essentially, this 
is the co-action $C^\infty(M)
\to C^\infty(M)\otimes C(\T^d)$ which yields the identification with 
$C^\infty(M)_\theta$.
Second, one uses the previous identification to deform the 
$C^\infty(M)$-module of 
sections $\Gamma(M,S)$ of the spinor bundle $S\to M$. 
Since the action $\tau$ is isometric, up to a double covering of 
$\T^d$,
the spinor bundle is $\T^d$-equivariant. The space of invariants 
$\Gamma(M,S)_\theta:=(\Gamma(M,S)\otimes C(\T^d)_\theta)^{\tau\times 
L^{-1}}$ is therefore 
a $C^\infty(M)_\theta$-module stable by the operator $D\otimes I$.
The restriction of the latter on 
$\Gamma(M,S)_\theta$
then defines the deformed Dirac operator $D_\theta$.

\vspace{1cm}

We end this subsection by stressing the fact that the way Connes and 
Dubois-Violette
define a deformation of the Dirac operator via an action of $\T^d$ 
essentially relies 
on two crucial properties. First, the invariance of Rieffel's 
deformed product
on $C(\T^d)_\theta$ under the (left or right) regular representation 
(this allows
to define the algebra (or module) structure on the space of invariant 
elements
in the tensor product). Second, the fact that the action $\tau$ of 
$\T^d$ on $M$ is isometric.

\subsection{Symmetric spaces and universal deformation formulae}
In this section we recall a consruction of a universal (strict) deformation 
formula for the actions of a non-Abelian Lie group. This formula was 
obtained via geometric methods based on symplectic symmetric spaces.
We begin by recalling the notion of symplectic symmetric space, we 
then pass to the specific example which will be relevant here.

A {\bf symplectic symmetric space} \cite{B95,BCG95} is triple $(M,\omega,s)$ where 
$(M,\omega)$ a symplectic manifold and $s:M\times M\to M$ is a 
smooth map such that for each point $x\in M$ the partial map $s_{x}:M\to M:y\mapsto 
s_x(y):=s(x,y)$ is an involutive  symplectic diffeomorphism which admits $x$ as 
isolated fixed point.
Under these hypotheses, the following formula 
defines a symplectic torsionfree affine connection 
$\nabla$ on $M$:
\begin{equation}
\omega_x(\nabla_XY,Z):=\frac{1}{2}X_x.\omega(Y+
s_{x\star}Y,Z);
\end{equation}
where $X,Y$ and $Z$ are smooth vector fields on $M$.
This connection is the only one being invariant 
under the {\sl symmetries} $\{s_x\}_{x\in M}$ 
which turns out to be the geodesic symmetries.
The group $G(M)$ generated by the products $s_x\circ
s_y$ is a Lie group of
 transformations of $M$ which acts transitively 
 on $M$. It is called the {\bf transvection } 
 group of $M$.
 
The particular example we will be concerned with  is the one where the transvection group is the Poincar\'e group $P_{1,1}:=
SO(1,1)\times\R^2$. It is a three dimensional solvable Lie group whose generic (two dimensional)
coadjoint orbits are symplectic symmetric spaces admitting $P_{1,1}$ as tranvection group.
They are  all equivalent to oneanother ($P_{1,1}/\R$ as homogeneous spaces) and topologically $\R^2$.

Let $(M,\omega)$ be such an orbit and let $\nabla$ be tis canonical connection. 
Such a space turns out to be stricly geodesically convex in the sense that between two points in $M$
there is a unique geodesic. It has moreover the property that given any three points $x,y$ and $z$ in $M$,
there's a unique fourth point $t$ in $M$ satisfying
\begin{equation}
s_xs_ys_zt=t.
\end{equation}
This allows to define the {\bf three point phase} of $M$ as the three point function $S\in C^\infty(M\times M\times M,\R)$ \cite{W94}
given by
\begin{equation}
S(x,y,z):=\int_{\Delta(t,s_zt,s_ys_zt)}\omega
\end{equation}
where $\Delta(x,y,z)$ denotes the geodesic triangle with vertices $x,y$ and $z$. The three point {\bf amplitude}
$A\in C^\infty(M\times M\times M,\R)$ \cite{Q97,BDRS04} is defined as the ratio
\begin{equation}
A(x,y,z):=\left(\frac{\int_{\Delta(t,s_zt,s_ys_zt)}\omega}{\int_{\Delta(x,y,z)}\omega}\right)^{\frac{1}{2}}.
\end{equation}
Both $S$ and $A$ are invariant under the diagonal action of the symmetries. For compactly  supported $u,v\in C^\infty_c(M)$,
one sets
\begin{equation}
u\star_\theta v(x):=\frac{1}{\theta^2}\int_{M\times M}A(x,y,z)\,e^{\frac{i}{\theta}S(x,y,z)}u(y)\,v(z)\,dy\,dz,
\mbox{ or shortly }u\star_\theta v=\frac{1}{\theta^2}\int_{M\times M}A\,e^{\frac{i}{\theta}S}\,u\otimes v;
\end{equation}
where $dx$ denotes the symplectic measure on $(M,\omega)$. This formula defines a symmetry invariant associative (non-formal) 
deformation
of $M$ led by the Poisson structure associated to $\omega$ on some topological function space ${\cal E}_\theta(M)\subset
C(M)$ \cite{}. For every real value of $\theta$ the space ${\cal E}_\theta(M)$ 
contains the test functions. Moreover, the deformed product is {\sl 
strongly closed} in the sense that, one has the trace property:
\begin{equation}\label{TRACE}
\int_{M}u\star_\theta v=\int_M uv
\end{equation}
for all compactly supported smooth functions $u$ and $v$.

We now explain how this is relevant to our situation of BHTZ spaces. The transvection group $P_{1,1}$ contains a copy of the group 
$\RR$ which acts
on maximal extensions of BHTZ spaces as discussed above. Moreover the 
restricted action $\RR\times M\to M$
turns out to be strictly transitive, so that one has the identification 
$M\simeq\RR$. Transporting the above structure to the group manifold $\RR$,
one gets a left-invariant deformation, again denoted by 
$\star_{\theta}$, on the symplectic group $(\RR,\omega^\RR)$. Let us 
denote by $K^{\RR}_\theta$ the associated (oscillating)  kernel:
\begin{equation}
u\star_\theta v=\int_{\RR\times\RR}K^\RR_\theta \; u\otimes v\quad u,v\in{\cal E}_\theta(\RR),
\end{equation}
where the integration is taken w.r.t. a left invariant  Haar measure on $\RR$.
Now, in a similar spirit as Rieffel's machinery, this defines a universal
deformation formula for actions of $\RR$. Indeed, if ${\cal X}$ is a space the group $\RR$
acts on via an action $\tau$, one defines for $a\in C({\cal X})$ and $x\in {\cal X}$ the function 
$\alpha^x:\RR\to\C:g\mapsto a(\tau_{g^{-1}}x)$. Setting 
\begin{equation}
{\cal E}_\theta({\cal X}):=\{a\in C({\cal X}):\alpha^xa\in
{\cal E}_\theta(\RR)\;\forall x\in {\cal X}\},
\end{equation}
 the following formula yields an associative product on the space
${\cal E}_\theta(\RR)$:
\begin{equation}
(a\star_\theta b)(x):=\int_{\RR\times\RR}K^\RR_\theta(e,g,h)\,\alpha^xa(g)\,\alpha^xb(h)\,dg\,dh;
\end{equation}
or shortly: $a\star_\theta b=:\int_{\RR\times\RR}K^\RR_\theta\;\alpha a\otimes\alpha b$.
Of course, such a product enjoys trace properties similar to (\ref{TRACE}).

Analoguous universal strict deformation formulae for more general 
solvable groups have been given in \cite{}.

\subsection{Spectral triples from non-Abelian Lie group actions}
We now give a noncommutative spectral triple deforming the commutative
one canonically associated to a maximally extended non-rotating 
BHTZ-space. We start, more generally,
with the data of a pseudo-Riemannian spin manifold $M$ of signature 
$(m,n)$. We set $D$ for the Dirac operator on $M$
acting on the $C^\infty(M)$-bi-module of smooth sections 
$\underline{S}$ of the spinor bundle $S\to M$ 
over $M$\footnote{More generally,  given a vector bundle $E\to M$ 
over $M$, we shall denote its space of smooth sections by 
$\underline{E}$.}. We note that the fact that, for all $a$ in 
$C^\infty(M)$, $[D,a]$ is continuous with respect 
to any reasonable topology on the spinor bi-module relies on the 
property that 
$D$ is a derivation of the module structure. In the case of Connes 
and 
Landi's construction, this property of the (undeformed) Dirac 
operator remains 
through the deformation because
the torus is an Abelian group. One difficulty in the case of a 
non-Abelian 
Lie group action is therefore to force the derivation property of the 
Dirac 
operator with respect to the deformed module structure. 

\subsubsection{Deformations of the spinor module and module 
endomorphisms}
Assume our deformation at the level of functions on $M$ comes from an 
action $\alpha$ of a symplectic solvable 
group $\RR$ in such a way that the spinor module is $\RR$-equivariant.
This means that there exists an action 
$\alpha^S:\RR\times\underline{S}\to\underline{S}$
such that $\alpha_g(a.\Psi)=\alpha_g(a).\alpha_g^S(\Psi)$ for all
$a\in C^\infty(M)$. Denoting by $K_\theta$ the convolution kernel 
which defines the deformation, 
the following formula yields, at least formally, a deformation of the 
(right) module structure:
\begin{equation}
\Psi\star_\theta a:=\int_{\RR\times\RR}\,K_{\theta}(e,g,h)\,
\alpha_h(a).\alpha^S_g(\Psi)\,\mbox{d}g\,\mbox{d}h\,;
\end{equation}
where $e$ denotes the identity in $\RR$. One writes a similar formula 
$a\star_\theta\Psi $ for a deformation of the left module structure. 
In the same way, one deforms the module structure attached to any 
vector bundle associated
to the spin structure over $M$. Note that, for the particular case of 
the endomorphism bundle,
one formally has, for all $\Psi\in\underline{S}, 
\,\gamma\in\underline{S^\star\otimes S}$ and $a\in C^\infty{M}$,
\begin{equation}
(\gamma\star_\theta \Psi )\star_\theta a = \gamma\star_\theta(\Psi\star_\theta 
a).
\end{equation}
In other words, the left deformed module of endomorphisms acts on 
the right deformed spinor module.

\subsubsection{Deforming the Dirac operator: case of trivial fibrations}
First consider the case of a differential operator 
$\Theta:\underline{S}\to\underline{S}$ of the form
\begin{equation}
\Theta(\Psi) = \gamma\left(\nabla_{X}\Psi\right),
\end{equation}
where $\nabla$ denotes the spin connection, where $X$ is a vector 
field on $M$ and where $\gamma\in\underline{S^\star\otimes S}$ is
a section of the endomorphism bundle.
Note that $\Theta$ acts as a derivation
of the (undeformed) module structure. \\
If the vector field $X$ commute with the actions $\alpha$ and 
$\alpha^S$, i.e. $X.\alpha_ga=\alpha_g(X.a)$ and 
$\nabla_X\alpha^S_g\Psi=\alpha^S_g\left(\nabla_X\Psi\right)$, then 
it automatically acts as derivation of the deformed multiplications. 
Indeed,
one formally has: $X.(a\star_\theta b)=\int K_\theta\,X.(\alpha 
a\otimes\alpha b)=
\int K_\theta\,(\alpha(X.a)\otimes \alpha b +\alpha 
a\otimes\alpha(X.b))=(X.a)\star_\theta b+a\star_\theta(X.b)$,
and similarly for $\nabla_X(\Psi\star_\theta a)$. In particular, in 
the case the element $\gamma$ is constant,
the operator $\Theta$ is then a derivation of the deformed module 
structure.  

Now, we turn back to our (extended) BHTZ spaces.
In the preceeding sections, we have seen that their maximal extension domains
are total spaces of trivial principal fibrations over $\R$ with structure group $\RR$.
In particular, this implies that the space ${\cal E}_\theta({\cal U})$, resulting from the 
application of the deformation formula to the left action $\tau$ of $\RR$ by isometries 
on ${\cal U}$, is non-trivial. Moreover, the (non-isometric) {\sl right} action of $\RR$
on itself--- and hence on a maximal extension domain ${\cal U}$--- defines 
for any element $X$ in the Lie algebra of the group $\RR$ a fundamental vector field again denoted by $X$ which 
commutes with the actions $\alpha$ and $\alpha^S$ of $\RR$. Observing that
$\tau_\star X=X$, one immediately concludes in the case of the non-rotating case $J=0$
\footnote{In the rotating case, the lack of a transversal everywhere orthogonal to the foliation
makes the argument more involved.}
that the Dirac operator $D$ on ${\cal U}$ acts as a derivation of the (right) deformed module structure.
This yields the formula $[D,a].{\bf \Psi}={\bf \Psi}\star_\theta(Da)$, showing that for every $a$ the commutator $[D,a]$
is a bounded operator in the space ${\cal H}$ of $L^2$ spinors.
The triple $({\cal E}_\theta({\cal U}), {\cal H}, D)$ is thus a spectral triple in the sense of Connes\footnote{More precisely,
because of noncompactess, in the sense of \cite{}.}. A distance formula for such Lorentzian spectral triples
has been formulated by V. Moretti in \cite{}. In a future work, we plan to investigate the question of understanding
how causality is encoded at the operator algebraic level.
\begin{rmk}
{\rm 
\begin{enumerate}
\item[(i)] In the more general case of a non-constant
element $\gamma$, the deformed operator $\Theta_\theta$ defined as 
\begin{equation}
\Theta_\theta(\Psi):=\left(\nabla_{X}\Psi\right)\star_\theta\gamma,
\end{equation}
acts as a derivation of the left deformed module structure (as a 
immediate consequence of the fact that the right deformed module of 
endomorphisms acts on the left deformed spinor module).\\
\item[(ii)] Another situation of interest is when the vector field $X$ is 
Hamiltonian with respect to the Poisson structure induced
by the action $\alpha$. Assume moreover that the spinor bundle $S\to 
M$ is trivial---as it is the case for $\ads$ or BHTZ-spaces.
In this case, one has a globally defined connection form 
$\Gamma\in\Omega^1(S^\star\otimes S)$ such that
\begin{equation}
\nabla_X\Psi=X.\Psi +\Gamma(X)\Psi.
\end{equation}
Denoting by $\lambda_X\in C^\infty(M)$ the hamiltonian function 
associated to $X$, one 
deforms the covariant derivative in the direction of $X$ as
\begin{equation}
(\nabla_X)^\theta\Psi:=\lambda_X\star_\theta\Psi-\Psi\star_\theta\lambda_X+\Psi\star_\theta\Gamma(X).
\end{equation}
Note that $(\nabla_X)^\theta$ provides a derivation of the left 
deformed module structure in the sense that  
$(\nabla_X)^\theta(a\star_\theta\Psi)=
[\lambda_X,a]_\theta\star_\theta\Psi+a\star_\theta(\nabla_X)^\theta\Psi$ 
where $[a,b]_\theta:=a\star_\theta b-b\star_\theta a$.
\end{enumerate}
}
\end{rmk}
These remarks have been used in \cite{} to define strict deformations of BHTZ spaces themselves 
rather than their extension domains.

\end{document}